\documentclass[10pt]{article}
\usepackage{amsmath}
\usepackage{amsfonts}
\usepackage{amssymb}
\usepackage[a 4 paper, left=1 in, right=1 in, top=1 in, foot skip=0.5 in ]{geometry}
\begin{document}
\setlength{\baselineskip}{0.7 cm}
\begin{Large}
\begin{center}
\textbf{On $\acute{C}$iri$\acute{c}$-type Contractive Mappings and Fixed Points over a $C^*-$algebra Valued Metric Space}
\end{center}
\end{Large}
\textbf{Kushal Roy$^1$ . Mantu Saha$^1$}\\\\
\textbf{Abstract}  In this paper $\acute{C}$iri$\acute{c}$ contractive mappings of type1 and type2 over a $C^*-$algebra valued metric space have been studied with illustrative examples. Some fixed point theorems with necessary and sufficient conditions have been proved for $\acute{C}$iri$\acute{c}$-type contractive mappings in a $C^*-$algebra valued metric space.\\\\
\textbf{Keywords}  $C^*-$algebra . $C^*-$algebra valued metric space . $\acute{C}$iri$\acute{c}$-type1 contractive mapping . $\acute{C}$iri$\acute{c}$-type2 contractive mapping . orbitally continuous . fixed point\\\\  
\textbf{AMS Mathematics Subject Classification(2010)} 47H10 . 54H25\\\\\\\\\\\\\\\\\\\\
\textbf{Acknowledgements} The first author wishes to thank the Council of Scientific and Industrial Research, New Delhi, India, for providing the research fellowship to continue the research work.\\
-----------------------------\\
Kushal Roy\\
Email: kushal.roy93@gmail.com\\\\
Mantu Saha\\
Email: mantusaha.bu@gmail.com\\\\
1 Department of Mathematics, The University of Burdwan, Burdwan-713104, West Bengal, India\\\\
\textbf{1 Introduction}\\
The most celebrated work 'Banach Contraction Principle'(See [2]) in proving the existence of a unique fixed point for a contraction mapping over a complete metric space with it's application has been established in many directions, namely (i) the usual contractive condition is relaxed by weakly contractive condition and (ii) the underlying spaces are replaced by metric spaces endowed with ordered or partially ordered structure. Recently O'Regan and Petrusel[11] had started their investigations towards the existence of a fixed point for generalized contraction mappings over an ordered metric space. In recent past years Caballero et al.[4] had been succeeded to apply his fixed point theorems particularly in a ordered metric space, for the existence of a solution of an ordinary differential equation. Subsequently with the introduction of cone metric space by Huang and Zhang[6] many researchers took an attempt to prove the generalized fixed point theorems on different cone metric spaces(See [1],[7],[12]). In 1971 L. B. $\acute{C}$iri$\acute{c}$ had defined the contraction type mapping(See [5]). It is well known that in a metric space a contraction mapping is always a $\acute{C}$iri$\acute{c}$ contraction type mapping but the converse is not true in general. Influenced on the survey of literatures of $C^*-$algebra(See [10]) and $C^*-$algebra valued metric space(see [9]) some fixed point theorems for $\acute{C}$iri$\acute{c}$-type contractive mappings on a $C^*-$algebra valued metric space have been proved by us. In this paper we also prove some common fixed point theorems for two mappings over such a $C^*-$algebra valued metric space. Strength of the hypothesis made in theorems has been examined by illustrative examples.\\ Before going to our main results we now need the following basic preliminary ideas. In the year 2014 Ma et al.[9] introduced the concept of $C^*-$algebra valued metric space. Here we denote an unital $C^{*}$-algebra by $\mathbb{A} $. An element  $a  \in \mathbb{A} $ is said to be a positive element if $a^*=a$ and $ \sigma (a) \subset \mathbb{R}^{+} $, where $ \sigma (a)= \lbrace \lambda \in \mathbb{C}: \lambda I-a$ is not invertible, $I$ is the unity in $ \mathbb{A}  \rbrace $. If $a$ is positive we write  it as $a\succeq \theta $, where $\theta $ is the zero element of  $ \mathbb{A} $. A partial ordering on $ \mathbb{A} $ is defined by $a\succeq b$ if and only if $(a-b) \succeq \theta $ whenever $a,b \in \mathbb{A} $. We denote the set of all positive elements of $ \mathbb{A} $ by $ \mathbb{A}_+ $ that is $ \mathbb{A}_+ = \lbrace a \in \mathbb{A}: a\succeq \theta \rbrace $. Also by $\acute{\mathbb{A}_+} $ we denote the set of all positive elements in $ \mathbb{A} $ which commutes with each element of $ \mathbb{A} $. Each positive element of a $C^{*}$-algebra $ \mathbb{A} $ has a unique positive square root.\\\\
\textbf{2 Some Basic Preliminaries }\\
\textbf{Definition 2.1}[9] Let $X$ be a non-empty set. Suppose the mapping $ d:X \times X\rightarrow \mathbb{A} $ satisfying the following conditions.
\\
1. $d(x,y) \succeq \theta\hspace{.3 cm} \forall  x,y \in X$ and $d(x,y)=\theta$ if and only if $x=y$
\\
2. $d(x,y)=d(y,x)\hspace{.3 cm}\forall  x,y \in X$
\\
3. $d(x,y) \preceq d(x,z)+d(z,y) \hspace{.3 cm}\forall x,y,z \in X$ \\
Then $d$ is called a $C^{*}$-algebra valued metric on $X$ and $(X,\mathbb{A},d)$ is called a $C^{*}$-algebra valued metric space.\\\\
\textbf{Example 2.2}[9] Let $X=\mathbb{R}$ and $\mathbb{A}=M_{2\times 2}(\mathbb{R})$ with $||A||=max\{|a_1|,|a_2|,|a_3|,|a_4|\}$, where ${a_i}'$s are the entries of the matrix $A\in M_{2\times 2}(\mathbb{R})$. Then $(X,\mathbb{A},d)$ is a $C^*-$algebra valued metric space, where $d(x,y)=$ diag$(|x-y|,\alpha |x-y|)$ and the partial ordering on $\mathbb{A}$ is given by $$\left(\begin{matrix}
a_1 & a_2 \cr a_3 & a_4
\end{matrix}\right)\succeq  \left(\begin{matrix}
b_1 & b_2 \cr b_3 & b_4
\end{matrix}\right)$$ iff $a_i \geq b_i$ for $i=1,2,3,4.$
\\\\\textbf{Definition 2.3}[9]
A sequence $ \lbrace x_n\rbrace $ in a $C^{*}$-algebra valued metric space $(X,\mathbb{A},d)$ is called convergent and converges to some $x \in X$ with respect to $ \mathbb{A} $ if for any $\epsilon>0$ $\exists  N \in  \mathbb{N} $ such that
\begin{center}
$||d(x_n,x)||<\epsilon$ whenever $n\geq N$ i.e. $d( x_n ,x)\rightarrow \theta$ as $n\rightarrow \infty$ \end{center}and we write $x_n \rightarrow x$ as $n\rightarrow \infty$ in $X$.\\\\
A sequence $ \lbrace x_n \rbrace $ in a $C^{*}$-algebra valued metric space $(X,\mathbb{A},d)$ is said to be Cauchy with respect to $ \mathbb{A} $ if for any $\epsilon>0$ $ \exists  N \in $ $ \mathbb{N} $ such that
\begin{center}
$||d(x_n,x_m)||<\epsilon$ whenever $n, m\geq N$ i.e. $d( x_n,x_m )\rightarrow$ $\theta$ as $n, m\rightarrow$ $\infty$ 
\end{center}
A $C^{*}$-algebra valued metric space is called complete if any cauchy sequence in $X$ with respect to $ \mathbb{A} $ is convergent to an element of $X.$\\\\
\textbf{Definition 2.4}[3] Let $(X, \mathbb{A} ,d)$ be a $C^{*}$-algebra valued metric space. A mapping $T:X\rightarrow X$ is said to be continuous at $a\in X$ with respect to $\mathbb{A}$ if for any $\epsilon>0$ $\exists$ $\delta>0$ such that $||d(Tx,Ta)||<\epsilon$ whenever $||d(x,a)||<\delta$. T is called continuous on $X$ with respect to $\mathbb{A}$ if it is continuous at any point $a\in X$. Evidently if $\lbrace x_n \rbrace$ be a sequence in $X$ such that $\lbrace x_n \rbrace \rightarrow x$ as n$\rightarrow \infty $ $\Rightarrow$ $T x_n \rightarrow Tx$ as $n\rightarrow \infty$ in $X$ that is $d(x_n,x)\rightarrow$ $\theta$ $ \Rightarrow $ $d(T x_n,Tx)\rightarrow$ $\theta$ as $n\rightarrow \infty$ then $T$ is continuous at $x \in X$.\\\\
\textbf{Lemma 2.5}[8]
We state some results that hold in a $C^{*}$-algebra with the partial ordering generated by the positive elements\\
1. $\mathbb{A}_+$=$\lbrace a^{*}a: a\in \mathbb{A} \rbrace$\\
2. If $a^{*}=a$, $b^{*}=b$, $a\preceq b$ and $c\in$ $\mathbb{A}$ then $c^{*}ac$ $\preceq$ $c^{*}bc$\\
3. For all $a,b\in$ $\mathbb{A}$ such that $a^{*}=a$, $b^{*}=b$ if $\theta$ $\preceq a\preceq b$ then $||a|| \leq ||b||$\\\\
\textbf{Lemma 2.6}[9]
Suppose that $\mathbb{A}$ is a unital $C^{*}$-algebra with a unity $I$.
\\1. If $a\in$ $\mathbb{A}_+$ with $||a||<\frac{1}{2}$, then $I-a$ is invertible and $||a(I-a)^{-1}||<1$\\
2. Suppose that $a,b\in$ $\mathbb{A}$ with $a,b\succeq$ $\theta$ and $ab=ba$ then $ab\succeq$ $\theta$\\
3. By $\acute{\mathbb{A}}$ we denote the set $\lbrace a \in \mathbb{A}: ab=ba \hspace{.2 cm} \forall b \in \mathbb{A} \rbrace$. Let $a\in$ $\acute{\mathbb{A}}$. If $b,c\in$ $\mathbb{A}$ with $b\succeq c\succeq$ $\theta$ and $I-a\in$ $\acute{\mathbb{A}}_+$ is invertible, then $(I-a)^{-1}b$ $\succeq$ $(I-a)^{-1}c$
 \\\\\textbf{Lemma 2.7}[10] Let $\mathbb{A}$ be an unital $C^*$-algebra with unity $I$\\ If $a,b$ are positive elements of $\mathbb{A}$, then the inequality $a\preceq b$ implies $a^\frac{1}{2}\preceq b^\frac{1}{2}$.\\\\
\textbf{Definition 2.8}[9] If $(X,\mathbb{A},d)$ is a $C^*-$algebra valued metric space and $T:X\rightarrow X$ satisfies\begin{equation}
d(Tx,Ty)\preceq a^* d(x,y) a\hspace{.2 cm},a\in \mathbb{A},\hspace{.2 cm}||a||<1\hspace{.2 cm}; \forall x,y\in X
\end{equation} then $T$ is called a $C^*-$algebra valued Contractive mapping.\\\\
In the following section we prove some theorems in light of the work of $\acute{C}$iri$\acute{c}$(See [5]) and Saha and Baisnab(See [14], [15]).\\ 
\textbf{3 Main Results}\\
We begin this portion with the definition of contractive type mapping on a $C^*-$algebra valued metric space.\\\\
\textbf{Definition 3.1} Let $(X,\mathbb{A},d)$ be a $C^*-$algebra valued metric space. A mapping $T:X\rightarrow X$ is called a $C^*-$algebra valued $\acute{C}$iri$\acute{c}$-type1 contractive mapping if for every $x,y\in X$ there are elements $q(x,y)$ in $\mathbb{A}$ satisfies $0\leq ||q(x,y)||< 1$ and there is a mapping $\delta:X\times X\rightarrow \mathbb{A}_+$ such that \begin{equation}
d(T^n x,T^n y)\preceq (q(x,y)^*)^n \delta(x,y) q(x,y)^n
\end{equation} holds for every $n\in \mathbb{N}.$ \\
Let us denote the collection of all $C^*-$algebra valued $\acute{C}$iri$\acute{c}$-type1 contractive mappings over $(X,\mathbb{A},d)$ by ${Ci}^*_1(X).$\\\\
\textbf{Definition 3.2} Let $(X,\mathbb{A},d)$ be a $C^*-$algebra valued metric space. A mapping $T:X\rightarrow X$ is called a $C^*-$algebra valued $\acute{C}$iri$\acute{c}$-type2 contractive mapping if for every $x,y\in X$ there are elements $q(x,y)$ in $\acute{\mathbb{A}_+}$, $0\leq ||q(x,y)||< 1$ and there is a mapping $\delta:X\times X\rightarrow \mathbb{A}_+$ such that \begin{equation}
d(T^n x,T^n y)\preceq q(x,y)^n \delta(x,y)
\end{equation} holds for every $n\in \mathbb{N}.$ \\
Let us denote the collection of all $C^*-$algebra valued $\acute{C}$iri$\acute{c}$-type2 contractive mappings over $(X,\mathbb{A},d)$ by ${Ci}^*_2(X).$\\\\
Note that in an abelian $C^*-$algebra valued metric space $(X,\mathbb{A},d)$ a $C^*-$algebra valued $\acute{C}$iri$\acute{c}$-type1 contractive mapping is a $C^*-$algebra valued $\acute{C}$iri$\acute{c}$-type2 contractive mapping and vice-versa. So for a commutative $C^*-$algebra valued metric space $(X,\mathbb{A},d)$ we have ${Ci}^*_1(X)={Ci}^*_2(X).$\\\\
\textbf{Definition 3.3}[13] Let $(X,\mathbb{A},d)$ be a $C^*-$algebra valued metric space and $T$ be a self map on $X.$ Then $T$ is said to be orbitally continuous at $u\in X$ if for any $x\in X$ $||d(T^{n_i}x,u)||\rightarrow 0$ as $i\rightarrow \infty$ implies $||d(T^{n_i +1}x,Tu)||\rightarrow 0$ as $i\rightarrow \infty$.\\\\
\textbf{Example 3.4} Let $X=[0,1]$ and $\mathbb{A}=M_2(\mathbb{R})$, where $$\left(\begin{matrix}
a & b\cr c & d
\end{matrix}\right)^*=\left(\begin{matrix}
a & b\cr c & d
\end{matrix}\right)^t $$, with $||A||=max\{ |a_1|,|a_2|,|a_3|,|a_4|\}$ where ${a_i}$'s are the entries of the matrix $A\in M_2(\mathbb{R})$. We define a partial ordering on $\mathbb{A}$ by $$A=\left(\begin{matrix}
a_1 & a_2 \cr a_3 & a_4
\end{matrix}\right)\succeq  \left(\begin{matrix}
b_1 & b_2 \cr b_3 & b_4
\end{matrix}\right)=B$$ iff $a_i \geq b_i$ for $i=1,2,3,4$ and $d:X\times X\rightarrow \mathbb{A}_+$ is defined by $$d(x,y)=\left(\begin{matrix}
|x-y| & 0 \cr 0 & |x-y|
\end{matrix}\right), \forall x,y\in X.$$ Then $(X,\mathbb{A},d)$ is a $C^*-$algebra valued metric space.\\Let $T:X\rightarrow X$ be defined by $Tx=\frac{x}{5}\hspace{.2 cm}\forall x\in X,$\\$\delta:X\times X\rightarrow \mathbb{A}_+$ be defined by $$\delta(x,y)=\left(\begin{matrix}
x+y & 0\cr 0 & x+y
\end{matrix}\right)\forall x,y\in X$$ and\\ $q:X\times X\rightarrow \mathbb{A}$
be defined by $$q(x,y)=\left(\begin{matrix}
\frac{1}{\sqrt{5}} & 0\cr 0 & \frac{1}{\sqrt{5}}
\end{matrix}\right)\forall x,y\in X$$ 
Then \begin{eqnarray} d(T^n x,T^n y)=\left(\begin{matrix}
|T^n x-T^n y| & 0\cr 0 & |T^n x-T^n y|
\end{matrix}\right)=\left(\begin{matrix}
\frac{1}{5^n}|x-y| & 0\cr 0 & \frac{1}{5^n}|x-y|
\end{matrix}\right)\nonumber\\=\left(\begin{matrix}
\frac{1}{\sqrt{5}} & 0\cr 0 & \frac{1}{\sqrt{5}}
\end{matrix}\right)^n\left(\begin{matrix}
|x-y| & 0\cr 0 & |x-y|
\end{matrix}\right)\left(\begin{matrix}
\frac{1}{\sqrt{5}} & 0\cr 0 & \frac{1}{\sqrt{5}}
\end{matrix}\right)^n\nonumber\\\preceq \left(\begin{matrix}
\frac{1}{\sqrt{5}} & 0\cr 0 & \frac{1}{\sqrt{5}}
\end{matrix}\right)^n\left(\begin{matrix}
x+y & 0\cr 0 & x+y
\end{matrix}\right)\left(\begin{matrix}
\frac{1}{\sqrt{5}} & 0\cr 0 & \frac{1}{\sqrt{5}}
\end{matrix}\right)^n\nonumber\\=(q(x,y)^*)^n\delta(x,y)q(x,y)^n \hspace{.2 cm}\forall x,y\in X.\nonumber\end{eqnarray}
Therefore $T$ is a $C^*-$algebra valued $\acute{C}$iri$\acute{c}$-type1 contractive mapping on $(X,\mathbb{A},d)$.\\\\Example 3.4 also supports in favour of $C^*-$algebra valued $\acute{C}$iri$\acute{c}$-type2 contractive mapping over a $C^*-$algebra valued metric space.\\\\
There are some mappings which are $C^*-$algebra valued $\acute{C}$iri$\acute{c}$-type1 contractive mappings($C^*-$algebra valued $\acute{C}$iri$\acute{c}$-type2 contractive mappings) but not Contraction mapping over a $C^*-$algebra valued metric space.\\
\textbf{Example 3.5} Let $X=[0,1]$ be the metric space with usual metric. Then it is a $C^*-$algebra valued metric space with the $C^*-$algebra $\mathbb{A}=\mathbb{R}.$\\Let $T:X\rightarrow X$ be defined by $Tx=\begin{cases} 0\hspace{.2 cm} \forall x\in [0,1) \cr \frac{1}{2} \hspace{.2 cm}if\hspace{.2 cm} x=1\end{cases}$\\Then $T^n x=0$ $\forall n\geq 1$ and $\forall x\in [0,1)$ and $T1=\frac{1}{2}$, $T^n 1=0$ $\forall n\geq 2.$\\We define $q:X\times X\rightarrow \mathbb{A}$ by $q(x,y)=\frac{1}{\sqrt{2}}$ $\forall x,y\in X$ and $\delta:X\times X\rightarrow \mathbb{A}_+$ by $\delta(x,y)=x+y$ $\forall x,y\in X.$\\ Then for all $x,y\in [0,1)$\\$d(Tx,Ty)=0\leq \frac{1}{2}(x+y)=\frac{1}{\sqrt{2}}(x+y)\frac{1}{\sqrt{2}}=q(x,y)^*\delta(x,y)q(x,y)$\\If $x\in [0,1)$ and $y=1$ then\\$d(Tx,Ty)=d(0,\frac{1}{2})=\frac{1}{2}\leq \frac{1}{2}(x+1)=\frac{1}{\sqrt{2}}(x+1)\frac{1}{\sqrt{2}}=q(x,y)^*\delta(x,y)q(x,y)$\\Similar is the case for $x=1,$ $y\in [0,1).$\\If $x=1=y$ then \\$d(Tx,Ty)=0< \frac{1}{2}(1+1)=\frac{1}{\sqrt{2}}(1+1)\frac{1}{\sqrt{2}}=q(x,y)^*\delta(x,y)q(x,y)$\\Now for any $n\geq 2$ and $\forall x,y\in [0,1]$ we get\\$d(T^n x,T^n y)=0\leq {\frac{1}{2}}^n(x+y)={\frac{1}{\sqrt{2}}}^n(x+y){\frac{1}{\sqrt{2}}}^n=(q(x,y)^*)^n\delta(x,y)q(x,y)^n$\\Therefore $d(T^n x,T^n y)\preceq(q(x,y)^*)^n\delta(x,y)q(x,y)^n$ $\forall x,y\in [0,1]$ and $\forall n\in \mathbb{N}$.\\So $T$ is a $C^*-$algebra valued $\acute{C}$iri$\acute{c}$-type1 contractive mapping($C^*-$algebra valued $\acute{C}$iri$\acute{c}$-type2 contractive mapping since the $C^*-$algebra $\mathbb{A}=\mathbb{R}$ is abelian) on $(X,\mathbb{A},d)$ but clearly it is not a contraction mapping.\\\\
We now go to our main theorems. Examples are also supported in support of our theorem.\\
\textbf{Theorem 3.6} Let $(X,\mathbb{A},d)$ be a complete $C^*-$algebra valued metric space, where $\mathbb{A}$ is a unital $C^*-$algebra. Let $T\in {Ci}^*_1(X)$. Then $T$ is orbitally continuous in $X$ if and only if $T$ has a unique fixed point in $X.$ \\\\
\textbf{Proof:} First let us assume that $T$ is orbitally continuous in $X.$ Let $x_0\in X$ be fixed.\\For any $n\in \mathbb{N}$ we get,\\$d(T^n x_0,T^{n+1} x_0)\preceq(q(x_0,T x_0)^*)^n\delta(x_0,T x_0)q(x_0,T x_0)^n$\\So for any $1\leq n<m$ we get,\\ $d(T^n x_0,T^m x_0)\preceq d(T^n x_0,T^{n+1} x_0)+d(T^{n+1} x_0,T^{n+2} x_0)+...+d(T^{m-1} x_0,T^m x_0)$
$$\begin{array}{lcll}
\Rightarrow d(T^n x_0,T^m x_0)&\preceq &(q(x_0,T x_0)^*)^n\delta(x_0,T x_0)q(x_0,T x_0)^n+ &\\ &&(q(x_0,T x_0)^*)^{n+1}\delta(x_0,T x_0)q(x_0,T x_0)^{n+1}+...&\\&&...+(q(x_0,T x_0)^*)^{m-1}\delta(x_0,T x_0)q(x_0,T x_0)^{m-1}&\\\\\Rightarrow ||d(T^n x_0,T^m x_0)||&\leq & ||q(x_0,T x_0)||^{2n} ||\delta(x_0,T x_0)||+...+||q(x_0,T x_0)||^{2m-2} ||\delta(x_0,T x_0)||&\\\\&=&||q(x_0,T x_0)||^{2n}\frac{1-||q(x_0,T x_0)||^{2m-2n}}{1-||q(x_0,T x_0)||^2}||\delta(x_0,T x_0)||&\\\\ &\leq &\frac{||q(x_0,T x_0)||^{2n}}{1-||q(x_0,T x_0)||^2}||\delta(x_0,T x_0)||
\end{array}$$
Since $0\leq ||q(x,y)||<1\hspace{.2 cm}\forall x,y\in X$ therefore $d(T^n x_0,T^m x_0)\rightarrow \theta$ as $n\rightarrow \infty$ and thus $\{T^n x_0\}$ is a Cauchy sequence in $X.$ Since $X$ is complete then $\{T^n x_0\}$ is convergent in $X$. Let $\{T^n x_0\}$ converges to $z\in X.$ Since $T$ is orbitally continuous on $X$ therefore we get $T^{n+1} x_0\rightarrow Tz$ as $n\rightarrow \infty$. \\Hence $Tz=z$ and $z$ is a fixed point of $T.$\\Now if possible let, $w$ be another fixed point of $T.$ Then $Tw=w$ and we get $T^n z=z$, $T^n w=w$ $\forall n\geq 1$\\Therefore $d(z,w)=d(T^n z,T^n w)\preceq (q(z,w)^*)^n\delta(z,w)q(z,w)^n\hspace{.3 cm}\forall n\in \mathbb{N}$\\Thus $||d(z,w)||\leq ||q(z,w)||^{2n}||\delta(z,w)||\rightarrow 0$ as $n\rightarrow \infty$. So $d(z,w)=\theta$ which implies $z=w.$\\Hence $T$ has a unique fixed point $z\in X.$ \\Conversely let, $T\in {Ci}^*_1(X)$ has a unique fixed point $u$ in $(X,\mathbb{A},d).$\\Let, $x_0\in X$ and $\{n_i\}$ be a subsequence of $\mathbb{N}.$\\Then, $d(T^{n_i} x_0, u)=d(T^{n_i} x_0,T^{n_i} u)\preceq (q(x_0,u)^*)^{n_i}\delta(x_0,u)q(x_0,u)^{n_i} $\\$\Rightarrow ||d(T^{n_i} x_0, u)||\leq ||q(x_0,u)||^{2 n_i} ||\delta(x_0,u)||\rightarrow 0$ as $i\rightarrow \infty$\\$\Rightarrow d(T^{n_i} x_0, u)\rightarrow \theta$ as $i\rightarrow \infty$. Therefore $lim_{i\rightarrow \infty}T^{n_i} x_0=u.$ \\Thus for any $x_0\in X$ and for any subsequence $\{n_i\}$ of $\mathbb{N}$ we get $lim_{i\rightarrow \infty}T^{n_i} x_0=u.$\\Now let, $m_i=n_i +1$ $\forall i\in \mathbb{N}$. Then $\{m_i\}$ is a subsequence of $\mathbb{N}$ and thus fom the above we get $lim_{i\rightarrow \infty}T^{m_i} x_0=u,$ implying that $lim_{i\rightarrow \infty}T T^{n_i} x_0=u=Tu.$\\So $T$ is orbitally continuous in $(X,\mathbb{A},d).$ This completes the proof.\\\\
\textbf{Theorem 3.7} Let $(X,\mathbb{A},d)$ be a complete $C^*-$algebra valued metric space, where $\mathbb{A}$ is a unital $C^*-$algebra. Let $T\in {Ci}^*_2(X)$. Then $T$ is orbitally continuous in $X$ if and only if $T$ has a unique fixed point in $X.$ \\\\
\textbf{proof:} The proof is similar to the Theorem 3.6.\\\\
\textbf{Note 3.8} We know that any continuous mapping on a $C^*-$algebra valued metric space $(X,\mathbb{A},d)$ is orbitally continuous, So a continuous mapping of ${Ci}^*_1(X)$ or ${Ci}^*_2(X)$ has a unique fixed point in a complete $C^*-$algebra valued metric space $X.$\\\\
\textbf{Note 3.9} Any orbitally continuous mapping may not be continuous in a $C^*-$algebra valued metric space $(X,\mathbb{A},d)$. Example 3.5 shows that $T$ is orbitally continuous in [0,1] but not continuous at $x=1.$ \\\\
\textbf{Example 3.10} Let $X=[0,1]$ be the metric space with usual metric. Then it is a complete $C^*-$algebra valued metric space with the $C^*-$algebra $\mathbb{A}=\mathbb{R}.$\\Let $T:X\rightarrow X$ be defined by $Tx=\frac{x}{3}$ $\forall x,y\in X$. Suppose $q:X\times X\rightarrow \mathbb{A}$ be defined by $q(x,y)=\frac{1}{\sqrt{3}}$ $\forall x,y\in X$ and $\delta:X\times X\rightarrow \mathbb{A}_+$ be defined by $\delta(x,y)=1+x+y$ $\forall x,y\in X$\\\\Then, $d(T^n x,T^n y)=|T^n x-T^n y|=\frac{1}{3^n}|x-y|\leq \frac{1}{3^n}(1+x+y)=(\frac{1}{\sqrt{3}})^n(1+x+y)(\frac{1}{\sqrt{3}})^n $\\$\Rightarrow d(T^n x, T^n y)\preceq (q(x,y)^*)^n\delta(x,y)q(x,y)^n$ $\forall x,y\in [0,1]$ and $\forall n\in \mathbb{N}$.\\Also clearly $T$ is orbitally continuous in $X$. Therefore all the conditions of theorem 3.6 are satisfied and we see that $0$ is the unique fixed point of $T$ in $X.$ \\\\
The condition that $T$ is orbitally continuous in Theorem 3.6 and Theorem 3.7 is necessary. which is evident from the following example\\\\
\textbf{Example 3.11} Let $X=\{0,\frac{1}{2},\frac{1}{2^2},\frac{1}{2^3}...\}$ be the metric space with usual metric. Then $X$ is a complete $C^*-$algebra valued metric space with the $C^*-$algebra $\mathbb{A}=\mathbb{R}.$\\\\Let $T:X\rightarrow X$ be defined by $Tx=\begin{cases}
\frac{1}{2}\hspace{.2 cm} if\hspace{.2 cm} x=0 \cr \frac{1}{2^{i+1}}\hspace{.2 cm} if\hspace{.2 cm}x=\frac{1}{2^i} \hspace{.2 cm}\forall i\geq 1.
\end{cases}$\\Define $q:X\times X\rightarrow \mathbb{A}$ by $q(x,y)=\frac{1}{\sqrt{2}}$ $\forall x,y\in X$ and $\delta:X\times X\rightarrow \mathbb{A}_+$ by $\delta(x,y)=1+x+y$ $\forall x,y\in X.$\\Let us choose $x_0=0$\\Then $T^n x_0=\frac{1}{2^n}$ $\forall n\in \mathbb{N}.$ So $T^n x_0\rightarrow 0$ as $n\rightarrow \infty$ and $T T^n x_0=T^{n+1} x_0\rightarrow 0\neq T0$.\\Therefore $T$ is not orbitally continuous in $X.$\\Now if $x=0$, $y=\frac{1}{2^k};k\geq 1$, then\\ $d(T^n x,T^n y)=d(\frac{1}{2^n},\frac{1}{2^{k+n}})=|\frac{1}{2^n}-\frac{1}{2^{k+n}}|=\frac{1}{2^n}|1-\frac{1}{2^k}|\leq (\frac{1}{\sqrt{2}})^n(1+x+y)(\frac{1}{\sqrt{2}})^n$\\$\Rightarrow d(T^n x, T^n y)\preceq (q(x,y)^*)^n\delta(x,y)q(x,y)^n$ $\forall n\in \mathbb{N}$.\\Also if $x=\frac{1}{2^i}$, $y=\frac{1}{2^j};\forall i,j\geq 1$, then\\ $d(T^n x,T^n y)=d(\frac{1}{2^{i+n}},\frac{1}{2^{j+n}})=|\frac{1}{2^{i+n}}-\frac{1}{2^{j+n}}|=\frac{1}{2^n}|\frac{1}{2^i}-\frac{1}{2^j}|\leq (\frac{1}{\sqrt{2}})^n(1+x+y)(\frac{1}{\sqrt{2}})^n$, which in turn implies that\\$d(T^n x, T^n y)\preceq (q(x,y)^*)^n\delta(x,y)q(x,y)^n$ $\forall n\in \mathbb{N}$. \\Therefore, $d(T^n x, T^n y)\preceq (q(x,y)^*)^n\delta(x,y)q(x,y)^n$ $\forall x,y\in X$ and $\forall n\in \mathbb{N}$.\\So, $T$ is a $C^*-$algebra valued $\acute{C}$iri$\acute{c}$-type1 contractive mapping($C^*-$algebra valued $\acute{C}$iri$\acute{c}$-type2 contractive mapping since the $C^*-$algebra $\mathbb{A}=\mathbb{R}$ is abelian) on $(X,\mathbb{A},d)$ but $T$ has no fixed point in $X.$\\\\
The condition that $X$ is complete in Theorem 3.6 and Theorem 3.7 is just sufficient.\\
\textbf{Example 3.12} Let $X=(-1,1)$ be the metric space with usual metric. Then $X$ is a $C^*-$algebra valued metric space with the $C^*-$algebra $\mathbb{A}=\mathbb{R}.$\\ Let $T:X\rightarrow X$ be defined by $Tx=\frac{x}{2}$ $\forall x\in X$. Define $q:X\times X\rightarrow \mathbb{A}$ by $q(x,y)=\frac{1}{\sqrt{2}}$ $\forall x,y\in X$ and $\delta:X\times X\rightarrow \mathbb{A}_+$ by $\delta(x,y)=4+x+y$ $\forall x,y\in X.$\\Then, $d(T^n x,T^n y)=|T^n x-T^n y|=\frac{1}{2^n}|x-y|\leq \frac{1}{2^n}(4+x+y)=(\frac{1}{\sqrt{2}})^n (4+x+y)(\frac{1}{\sqrt{2}})^n$, implying that  $d(T^n x, T^n y)\preceq (q(x,y)^*)^n\delta(x,y)q(x,y)^n$ $\forall x,y\in X$ and $\forall n\in \mathbb{N}$.\\Clearly $T\in {Ci}^*_1(X)$(also $T\in {Ci}^*_2(X)$ since $\mathbb{A}=\mathbb{R}$ is abelian) and $T$ has a unique fixed point $0$ in $X$ but $X$ is not complete.\\\\
\textbf{Corollary 3.13} A $C^*-$algebra valued contractive mapping has always a unique fixed point in a complete $C^*-$algebra valued metric space.\\\\
\textbf{Proof:} Let $(X,\mathbb{A},d)$ be a $C^*-$algebra valued metric space. A mapping $T:X\rightarrow X$ is called a $C^*-$algebra valued Contractive mapping if it satisfies the condition $d(Tx,Ty)\preceq a^*d(x,y)a$ $\forall x,y\in X,$ where $a\in \mathbb{A}$ satisfies $0\leq ||a||<1.$\\ Then clearly $T$ is continuous on $X$ and thus orbitally continuous on $X$ and $d(T^n x,T^n y)\preceq (a^*)^n d(x,y)a^n$ $\forall x,y\in X$ and $\forall n\in \mathbb{N}.$\\Now if we define $q:X\times X\rightarrow \mathbb{A}$ by $q(x,y)=a$ $\forall x,y\in X$ and $\delta:X\times X\rightarrow \mathbb{A}_+$ by $\delta(x,y)=d(x,y)$ $\forall x,y\in X$ then we get $0\leq q(x,y)<1$ $\forall x,y\in X$ and $d(T^n x,T^n y)\preceq (q(x,y)^*)^n\delta(x,y)q(x,y)^n$ $\forall x,y\in X$ and $\forall n\in \mathbb{N}$. So $T\in {Ci}^*_1(X)$ and thus by Theorem 3.6 $T$ has a unique fixed point in $X.$\\\\
\textbf{Corollary 3.14} Let $(X,\mathbb{A},d)$ be a $C^*-$algebra valued metric space. Suppose the mapping $T:X\rightarrow X$ satisfies $\forall x,y\in X$ $d(Tx,Ty)\preceq A(d(x,Tx)+d(y,Ty))$, where $A\in \acute{\mathbb{A}_+}$ with $||A||<\frac{1}{2}.$ Then there exists a unique fixed point in $X.$\\\\
\textbf{Proof:} Let $x,y\in X$ be arbitrary. Then, $d(T^2 x,T^2 y)\preceq A(d(Tx,T^2 x)+d(Ty,T^2 y))$. Now, $d(Tx,T^2 x)\preceq A(d(x,Tx)+d(Tx,T^2 x))$. So $(I-A)d(Tx,T^2 x)\preceq A d(x,Tx)$. Hence $d(Tx,T^2 x)\preceq A(I-A)^{-1} d(x,Tx)=B d(x,Tx),$ where $B=A(I-A)^{-1}\in \acute{\mathbb{A}_+}$ and $||B||<1$ by Lemma 2.6. \\Similarly we get, $d(Ty,T^2 y)\preceq B d(y,Ty)$. Therefore, $d(T^2 x,T^2 y)\preceq AB (d(x,Tx)+d(y,Ty))$. Similarly $d(T^3 x,T^3 y)\preceq AB (d(Tx,T^2 x)+d(Ty,T^2 y))\preceq AB^2(d(x,Tx)+d(y,Ty))$\\Proceeding in this way we get, \\$d(T^n x,T^n y)\preceq AB^{n-1}(d(x,Tx)+d(y,Ty))\preceq B^n (d(x,Tx)+d(y,Ty))$ $\forall x,y\in X$ and $\forall n\in \mathbb{N}$ [$\because A\preceq B$]\\Therefore $T\in {Ci}^*_2(X)$\\Also if $u\in X$ and $lim_{i\rightarrow \infty}T^{n_i} x_0=u$ for some $x_0\in X$ then
$$\begin{array}{lcll}
d(T^{n_i +1}x_0,Tu)&\preceq & A(d(T^{n_i}x_0,T^{n_i +1}x_0)+d(u,Tu))&\\&\preceq & A(d(T^{n_i}x_0,T^{n_i +1}x_0)+d(u,T^{n_i +1}x_0)+d(T^{n_i +1}x_0,Tu))&\\&\preceq & A(I-A)^{-1}(d(T^{n_i}x_0,T^{n_i +1}x_0)+d(u,T^{n_i +1}x_0))&\\&\preceq& A(I-A)^{-1}(2 d(T^{n_i}x_0,T^{n_i +1}x_0)+d(u,T^{n_i}x_0))&\\&\preceq & A(I-A)^{-1}\{2 B^{n_i}(d(x_0,T x_0)+d(T x_0+T^2 x_0))+d(u,T^{n_i}x_0)\}&\\&\preceq & B\{2 B^{n_i}(I+B)d(x_0,T x_0)+d(u,T^{n_i}x_0)\}&\\
\end{array}$$
So $d(T^{n_i +1}x_0,Tu)\rightarrow \theta$ as $i\rightarrow \infty$.\\Therefore, $lim_{i\rightarrow \infty}T T^{n_i} x_0=Tu$. Hence $T$ is orbitally continuous in $X$ and thus by Theorem 3.7 $T$ has a unique fixed point in $X.$\\\\
\textbf{Theorem 3.15} Let $(X,\mathbb{A},d)$ be a complete $C^*-$algebra valued metric space, where $\mathbb{A}$ is a unital $C^*-$algebra. Two mappings $T,S:X\rightarrow X$ satisfy \\$d(T^n x,S^n y)\preceq (q(x,y)^*)^n\delta(x,y)q(x,y)^n$ $\forall x,y\in X$ and $\forall n\in \mathbb{N}$., where $q:X\times X\rightarrow \mathbb{A}$ be a mapping which satisfy $0\leq ||q(x,y)||<1$ $\forall x,y\in X$ and $\delta:X\times X\rightarrow \mathbb{A}_+$ be another mapping. Then $T$ and $S$ both are orbitally continuous in $X$ if and only if $T$ and $S$ have a unique common fixed point in $X.$\\\\
\textbf{Proof:} First let us assume that both $T$ and $S$ be orbitally continuous.\\Let $x_0\in X$ be fixed.\\Let us consider a sequence $\{x_n\}$ in $X,$ where $x_n=\begin{cases}
T^n x_0\hspace{.3 cm} $if $n$ is even$\cr S^n x_0 \hspace{.3 cm} $if $n$ is odd$ 
\end{cases}$ \\\\Then, $d(x_n,x_{n+1})=d(T^n x_0,S^{n+1} x_0)\preceq (q(x_0,S x_0)^*)^n\delta(x_0,S x_0)q(x_0,S x_0)^n$, whenever $n$ is even\\ and $d(x_n,x_{n+1})=d(S^n x_0,T^{n+1} x_0)=d(T^{n+1} x_0,S^n x_0)\preceq (q(T x_0, x_0)^*)^n\delta(T x_0, x_0)q(T x_0, x_0)^n$, whenever $n$ is odd.\\Therefore, $d(x_n,x_{n+1})\preceq (q(x_0,S x_0)^*)^n\delta(x_0,S x_0)q(x_0,S x_0)^n +(q(T x_0, x_0)^*)^n\delta(T x_0, x_0)q(T x_0, x_0)^n$ $\forall n\in \mathbb{N}$.\\Hence $||d(x_n,x_{n+1})||\leq ||q(x_0,S x_0)||^{2n} ||\delta(x_0,S x_0)||+||q(T x_0, x_0)||^{2n} ||\delta(T x_0, x_0)||$.\\So, $||d(x_n,x_{n+1})||\rightarrow 0$ as $n\rightarrow \infty$ and thus $d(x_n,x_{n+1})\rightarrow \theta$ as $n\rightarrow \infty.$\\Now for $1\leq n<m$ we get,\\$d(x_n,x_m)\preceq d(x_n,x_{n+1})+d(x_{n+1},x_{n+2})+...+d(x_{m-1},x_m)\rightarrow \theta$ as $n\rightarrow \infty.$ Hence $\{x_n\}$ is Cauchy in $X.$ Since $(X,\mathbb{A},d)$ is complete then there exists a unique $z\in X$ such that $x_n\rightarrow z$ as $n\rightarrow \infty$. Therefore $\{T^{2n}x_0\}_{n\in \mathbb{N}}$ and $\{S^{2n-1}x_0\}_{n\in \mathbb{N}}$ both converge to $z.$\\Since $T$ and $S$ both are orbitally continuous in $X$ so $T^{2n+1}x_0\rightarrow Tz$ and $S^{2n} x_0\rightarrow Sz$ as $n\rightarrow \infty.$\\Therefore, $d(Tz,z)\preceq d(Tz,T^{2n+1}x_0)+d(T^{2n+1}x_0,S^{2n+1}x_0)+d(S^{2n+1}x_0,z)$\\$\Rightarrow d(Tz,z)\preceq d(Tz,T^{2n+1}x_0)+(q( x_0, x_0)^*)^{2n+1}\delta( x_0, x_0)q( x_0, x_0)^{2n+1}+d(S^{2n+1}x_0,z)$\\$\Rightarrow ||d(Tz,z)||\leq ||d(Tz,T^{2n+1}x_0)||+||q( x_0, x_0)||^{4n+2}||\delta( x_0, x_0)||+||d(S^{2n+1}x_0,z)||\rightarrow 0$ as $n\rightarrow \infty.$\\Thus, $||d(Tz,z)||=0\Rightarrow d(Tz,z)=\theta\Rightarrow Tz=z.$\\Similarly we can show that $Sz=z.$ Therefore $Tz=z=Sz$ and so $z$ is a common fixed point of $T$ and $S.$\\If possible let, $u$ be another common fixed point of $T$ and $S.$ Thus $Tu=u=Su.$\\Therefore, $d(u,z)=d(T^n u,S^n z)\preceq (q(u,z)^*)^n\delta(u,z)q(u,z)^n$ $\forall n\in \mathbb{N}$\\So, $||d(u,z)||\leq ||q(u,z)||^{2n}||\delta(u,z)||\rightarrow 0$ as $n\rightarrow \infty$. Thus, $||d(u,z)||=\theta\Rightarrow d(u,z)=\theta\Rightarrow u=z$\\Hence $T$ and $S$ have a unique common fixed point in $X.$\\Conversely let, $T$ and $S$ both have unique common fixed point in $X.$\\So there exists a unique $z\in X$ such that $Tz=z=Sz.$ We first show that, $z$ is the unique fixed point of $T.$\\Let, $w\in X$ be another fixed point of $T.$ Then $Tw=w.$\\So, $d(w,z)=d(T^n w,S^n z)\preceq (q(w,z)^*)^n\delta(w,z)q(w,z)^n$ $\forall n\in \mathbb{N}$. Hence $||d(w,z)||\leq ||q(w,z)||^{2n}||\delta(w,z)||\rightarrow 0$ as $n\rightarrow \infty.$ Therefore, $||d(w,z)||=0\Rightarrow d(w,z)=\theta\Rightarrow w=z.$ So, $z$ is the unique fixed point of $T.$ Similarly, $z$ is the unique fixed point of $S.$\\Let, $x_0\in X$ and $\{n_i\}$ be a subsequence of $\mathbb{N}.$\\Then, $d(T^{n_i} x_0, z)=d(T^{n_i} x_0,S^{n_i} z)\preceq (q(x_0,z)^*)^{n_i}\delta(x_0,z)q(x_0,z)^{n_i} $. So $||d(T^{n_i} x_0, z)||\leq ||q(x_0,z)||^{2 n_i} ||\delta(x_0,z)||\rightarrow 0$ as $i\rightarrow \infty$, implying that $d(T^{n_i} x_0, z)\rightarrow \theta$ as $i\rightarrow \infty$. Therefore $lim_{i\rightarrow \infty}T^{n_i} x_0=z.$ \\Thus for any $x_0\in X$ and for any subsequence $\{n_i\}$ of $\mathbb{N}$ we get $lim_{i\rightarrow \infty}T^{n_i} x_0=z.$\\Now let, $m_i=n_i +1$ $\forall i\in \mathbb{N}$. Then $\{m_i\}$ is a subsequence of $\mathbb{N}$ and thus fom the above we get $lim_{i\rightarrow \infty}T^{m_i} x_0=z.$ implying that  $lim_{i\rightarrow \infty}T T^{n_i} x_0=z=Tz.$ So $T$ is orbitally continuous in $(X,\mathbb{A},d).$\\Similarly we get, $S$ is also orbitally continuous in $X$ which proves our Theorem.\\\\
\textbf{Theorem 3.16} Let $(X,\mathbb{A},d)$ be a complete $C^*-$algebra valued metric space, where $\mathbb{A}$ is a unital $C^*-$algebra. Suppose that two mappings $T,S:X\rightarrow X$ satisfy \\$d(T^n x,S^n y)\preceq q(x,y)^n\delta(x,y)$ $\forall x,y\in X$ and $\forall n\in \mathbb{N}$., where $q:X\times X\rightarrow \acute{\mathbb{A}_+}$ be a mapping such that $0\leq ||q(x,y)||<1$ $\forall x,y\in X$ and $\delta:X\times X\rightarrow \mathbb{A}_+$ be another mapping. Then $T$ and $S$ both are orbitally continuous in $X$ if and only if $T$ and $S$ have a unique common fixed point in $X.$ \\\\
\textbf{Proof:} The proof is similar to the above theorem.\\\\
\textbf{Example 3.17} Let $X=[0,1]$ be the metric space with usual metric. Then it is a complete $C^*-$algebra valued metric space with the $C^*-$algebra $\mathbb{A}=\mathbb{R}.$ We take two mappings $T$ and $S$ on $X$ such that \\$T:X\rightarrow X$ is defined by $Tx=\frac{x}{2}$ and $S:X\rightarrow X$ is defined by $Sx=\frac{x}{4}$ $\forall x\in X.$\\ Also define $q:X\times X\rightarrow \mathbb{A}$ by $q(x,y)=\frac{1}{\sqrt{2}}$ $\forall x,y\in X$ and $\delta:X\times X\rightarrow \mathbb{A}_+$ by $\delta(x,y)=1+x+y$ $\forall x,y\in X.$\\Now, $d(T^n x,S^n y)=|T^n x-S^n y|=|\frac{x}{2^n}-\frac{y}{4^n}|=\frac{1}{2^n}|x-\frac{y}{2^n}|\leq \frac{1}{2^n}(1+x+y)=(\frac{1}{\sqrt{2}})^n(1+x+y)(\frac{1}{\sqrt{2}})^n$\\$\Rightarrow d(T^n x,S^n y)\preceq (q(x,y)^*)^n\delta(x,y)q(x,y)^n$ $\forall x,y\in X$, $\forall n\in \mathbb{N}$.\\Also both $T$ and $S$ are continuous in $X$. So $T$ and $S$ are orbitally continuous in $X.$ We see that $0$ is the unique common fixed point of $T$ and $S$ in $X.$\\\\
The condition that $X$ is complete in the Theorem 3.15 is just sufficient.\\
\textbf{Example 3.18} Let $X=(-1,1)$ be the metric space with usual metric. Then it is a $C^*-$algebra valued metric space with the $C^*-$algebra $\mathbb{A}=\mathbb{R}.$ We take two mappings $T$ and $S$ on $X$ such that \\$T:X\rightarrow X$ is defined by $Tx=\frac{x}{3}$ and $S:X\rightarrow X$ is defined by $Sx=\frac{x}{6}$ $\forall x\in X.$\\ Define $q:X\times X\rightarrow \mathbb{A}$ by $q(x,y)=\frac{1}{\sqrt{3}}$ $\forall x,y\in X$ and $\delta:X\times X\rightarrow \mathbb{A}_+$ by $\delta(x,y)=4+x+y$ $\forall x,y\in X.$\\Then, $d(T^n x,S^n y)=|T^n x-S^n y|=|\frac{x}{3^n}-\frac{y}{6^n}|=\frac{1}{3^n}|x-\frac{y}{2^n}|\leq \frac{1}{3^n}(4+x+y)=(\frac{1}{\sqrt{3}})^n(4+x+y)(\frac{1}{\sqrt{3}})^n$. So $d(T^n x,S^n y)\preceq (q(x,y)^*)^n\delta(x,y)q(x,y)^n$ $\forall x,y\in X$, $\forall n\in \mathbb{N}$.\\Though $(X,\mathbb{A},d)$ is incomplete, $T$ and $S$ have a unique common fixed point in $X$ which is $0.$\\\\
The condition that both $T$ and $S$ are orbitally continuous in $X$ is necessary in the Theorem 3.15\\
\textbf{Example 3.19} Let $X=\{0,\frac{1}{2},\frac{1}{2^2},\frac{1}{2^3}...\}$ be the metric space with usual metric. Then $X$ is a complete $C^*-$algebra valued metric space with the $C^*-$algebra $\mathbb{A}=\mathbb{R}.$\\Let $T:X\rightarrow X$ be defined by $Tx=\begin{cases}
0\hspace{.2 cm} if\hspace{.2 cm} x=0 \cr \frac{1}{2^{i+1}}\hspace{.2 cm} if\hspace{.2 cm}x=\frac{1}{2^i} \hspace{.2 cm}\forall i\geq 1.
\end{cases}$\\Also Let $S:X\rightarrow X$ be defined by $Sx=\begin{cases}
\frac{1}{2}\hspace{.2 cm} if\hspace{.2 cm} x=0 \cr \frac{1}{2^{i+1}}\hspace{.2 cm} if\hspace{.2 cm}x=\frac{1}{2^i} \hspace{.2 cm}\forall i\geq 1,
\end{cases}$ \\$q:X\times X\rightarrow \mathbb{A}$ by $q(x,y)=\frac{1}{\sqrt{2}}$ $\forall x,y\in X$ and $\delta:X\times X\rightarrow \mathbb{A}_+$ by $\delta(x,y)=1+x+y$ $\forall x,y\in X.$\\Now if $x=0$, $y=\frac{1}{2^k};k\geq 1$ \\Then, $d(T^n x,S^n y)=d(0,\frac{1}{2^{k+n}})=|0-\frac{1}{2^{k+n}}|=\frac{1}{2^{n+k}}\leq (\frac{1}{\sqrt{2}})^n(1+x+y)(\frac{1}{\sqrt{2}})^n$. Hence\\$d(T^n x, S^n y)\preceq (q(x,y)^*)^n\delta(x,y)q(x,y)^n$ $\forall n\in \mathbb{N}$.\\If $x=\frac{1}{2^k};k\geq 1$, $y=0$ then we get\\$d(T^n x,S^n y)=d(\frac{1}{2^{n+k}},\frac{1}{2^{n}})=|\frac{1}{2^n}-\frac{1}{2^{k+n}}|=\frac{1}{2^n}|1-\frac{1}{2^k}|\leq (\frac{1}{\sqrt{2}})^n(1+x+y)(\frac{1}{\sqrt{2}})^n$\\So $d(T^n x, S^n y)\preceq (q(x,y)^*)^n\delta(x,y)q(x,y)^n$ $\forall n\in \mathbb{N}$.  \\Also if $x=\frac{1}{2^i}$, $y=\frac{1}{2^j};\forall i,j\geq 1$\\Then, $d(T^n x,S^n y)=d(\frac{1}{2^{i+n}},\frac{1}{2^{j+n}})=|\frac{1}{2^{i+n}}-\frac{1}{2^{j+n}}|=\frac{1}{2^n}|\frac{1}{2^i}-\frac{1}{2^j}|\leq (\frac{1}{\sqrt{2}})^n(1+x+y)(\frac{1}{\sqrt{2}})^n$.\\ Thus we have $d(T^n x, S^n y)\preceq (q(x,y)^*)^n\delta(x,y)q(x,y)^n$ $\forall n\in \mathbb{N}$. \\Therefore, $d(T^n x, S^n y)\preceq (q(x,y)^*)^n\delta(x,y)q(x,y)^n$ $\forall x,y\in X$ and $\forall n\in \mathbb{N}$.\\Here $X$ is complete, $T$ is orbitally continuous but $S$ is not orbitally continuous. Here we see that $T$ and $S$ have no common fixed point in $X.$\\\\
\textbf{Corollary 3.20} Let $(X,\mathbb{A},d)$ be a complete $C^*-$algebra valued metric space, where $\mathbb{A}$ is a unital $C^*-$algebra. Let $T,S$ be two mappings on $X$ such that\\$d(T^n x,(ST)^n y)\preceq (q(x,y)^*)^n\delta(x,y)q(x,y)^n$ $\forall x,y\in X$ and $\forall n\in \mathbb{N}$., where $q:X\times X\rightarrow \mathbb{A}$ be a mapping which satisfy $0\leq ||q(x,y)||<1$ $\forall x,y\in X$ and $\delta:X\times X\rightarrow \mathbb{A}_+$ be another mapping. Assume that $T$ and $ST$ both are orbitally continuous in $X.$ Then $T$ and $S$ have a unique common fixed point in $X.$\\\\
\textbf{Proof:} From the Theorem 3.15 it follows that $T$ and $ST$ have a unique common fixed point in $X,$ say $z.$ Therefore, $Tz=z=(ST)z.$ Now, $(ST)z=z\Rightarrow S(Tz)=z\Rightarrow Sz=z.$ So, $Tz=Sz=z$. Thus $z$ is a common fixed point of $T$ and $S.$ Let $w$ be another common fixed point of $T$ and $S.$ Then $Tw=w=Sw.$ So we get, $(ST)w=S(Tw)=Sw=w$ Therefore, $Tw=(ST)w=w.$\\So $w$ is a common fixed point of $T$ and $ST$ in $X$. Since $S$ and $ST$ have a unique common fixed point in $X$ it follows that $z=w.$ So $z$ is the unique common fixed point of $T$ and $S$ in $X.$\\\\
\textbf{Corollary 3.21} Let $(X,\mathbb{A},d)$ be a complete $C^*-$algebra valued metric space, where $\mathbb{A}$ is a unital $C^*-$algebra. Let $T,S$ be two mappings on $X$ satisfying\\$d(T^n x,(ST)^n y)\preceq q(x,y)^n\delta(x,y)$ $\forall x,y\in X$ and $\forall n\in \mathbb{N}$., where $q:X\times X\rightarrow \acute{\mathbb{A}_+}$ be a mapping such that $0\leq ||q(x,y)||<1$ $\forall x,y\in X$ and $\delta:X\times X\rightarrow \mathbb{A}_+$ be another mapping. Assume that $T$ and $ST$ both are orbitally continuous in $X.$ Then $T$ and $S$ have a unique common fixed point in $X.$\\\\
\textbf{Proof:} The proof is similar to the above.\\\\
\begin{large}
\textbf{References}
\end{large}\\\\
1. Abbas, M., Jungck, G., Common fixed point results for noncommuting mapping with out continuity in cone metric space, \textit{J. Math. Anal. Appl. 341}, 416-420 (2008).\\2. Banach, S., Sur les operations dans les ensembles abstraits et leur application aux equations untegrales, \textit{Fund. Math., 3}, 133-181 (1922).\\3. Batul, S. and Kamran, T., $C^*-$valued contractive type mappings, \textit{Fixed Point Theory and Applications 2015}, 142, DOI 10.1186/s13663-015-0393-3 (2015).\\4. Caballero, J., Harjani, J., Sadarangani, K., Contractive-like mapping principles in ordered metric spaces and applications to ordinary differential equations, \textit{Fixed Point Theory and Applications 2010}, Article ID 916064(2010).\\5. $\acute{C}$iri$\acute{c}$, L.B., Generalized contractions and fixed point theorems, \textit{Publications De L'Institut Mathematique}, 26, 19-26 (1971).\\6. Huang, LG., Zhang, X., Cone metric spaces and fixed point theorems of contractive mappings, \textit{J. Math. Anal. Appl. 332(2)}, 1468-1476, (2007).\\7. Kadelburg, Z. and Radenovic, S., A note on various types of cones and fixed point results in cone metric  spaces, \textit{Asian Journal of Mathematics and Applications}, Article ID ama0104, 7 pages (2013).\\8. Kadelburg, Z. and Radenovic, S., Fixed point results in $C^*-$algebra-valued metric spaces are direct consequences of their standard metric counterparts, \textit{Fixed Point Theory and Applications 2016}, 53, DOI 10.1186/s13663-016-0544-1 (2016).\\9. Ma, Z., Jiang, L. and Sun, H., $C^*-$algebra-valued metric spaces and related fixed point theorems, \textit{Fixed Point Theory and Applications 2014}, 206, DOI 10.1186/1687-1812-2014-206 (2014).\\10. Murphy, G.J., $C^*-$algebras and Operator Theory, \textit{Academic Press}, Boston, (1990).\\11. O'Regan, D., Petrusel, A., Fixed point theorems for generalized contractions in ordered metric spaces, \textit{J. Math. Anal. Appl. 341(2)}, 1241-1252 (2008).\\12. Rezapour, S., Hamlbarani, R., Some notes on the paper 'Cone metric spaces and fixed point theorems of contractive mappings', \textit{J. Math. Anal. Appl. 345}, 719-724 (2008).\\13. Roy, K., Saha, M., Fixed point theorems for generalized contractive and expansive type mappings over a $C^*-$algebra valued metric space, (communicated).\\14. Saha, M., Baisnab, A.P., Fixed point of mappings with contractive iterate, \textit{Proc. Nat. Acad. Sci., India}, 63(A), IV, 645-650 (1993).\\15. Saha, M., Baisnab, A.P., On $\acute{C}$iri$\acute{c}$ contraction operator and fixed points, \textit{Indian Journal of Pure and Applied Mathematics, Springer}, 27(2), 177-182 (1996).

\end{document}